\documentclass{amsart}
\usepackage{amscd}
\newtheorem{thm}{Theorem}[section]
\newtheorem{lemma}[thm]{Lemma}
\newtheorem{prop}[thm]{Proposition}
\newtheorem{quest}[thm]{Question}
\newtheorem{cor}[thm]{Corollary}

\newtheorem{remark}[thm]{Remark}

\theoremstyle{definition}

\newtheorem{exam}[thm]{Example}

\numberwithin{equation}{section}

\newcommand{\chara}{\mathrm{char}\,}
\newcommand{\bight}{\mathrm{bight}\,}
\newcommand{\codim}{\mathrm{codim}\,}

\newcommand{\frm}{{\mathfrak m}}

\newcommand{\frp}{{\mathfrak p}}

\newcommand{\height}{\mathrm{height}\,}

\newcommand{\indeg}{\mathrm{indeg}\,}
\newcommand{\link}{\mathrm{link}}

\newcommand{\reg}{\mathrm{reg}\,}
\newcommand{\rt}{\mathrm{rt}\,}

\newcommand{\Star}{\mathrm{star}}

\newcommand{\Z}{\ensuremath{\mathbb Z}}

\newcommand{\bbZ}{\ensuremath{\mathbb Z}}

\newcommand{\bbinom}[2]{%
\genfrac{(}{)}{0pt}{}{#1}{#2}}

\newcommand{\sbinom}[2]{%
\text{\small $\displaystyle{\left(\!\!\!\begin{array}{c}#1 \\ #2 \end{array}\!\!\!
\right)}$}}

\begin{document}
\title{Stanley--Reisner rings with large multiplicities are Cohen--Macaulay}
\author{Naoki Terai}
\address{Department of Mathematics, Faculty of Culture and Education, 
         Saga University, Saga 840--8502, Japan} 
\email{terai@cc.saga-u.ac.jp}
\author{Ken-ichi Yoshida}
\address{Graduate School of Mathematics, Nagoya University,
         Nagoya  464--8602, Japan}
\email{yoshida@math.nagoya-u.ac.jp}
\date{Mar.10, 2005}
\begin{abstract}
We prove that certain class of 
Stanley--Reisner rings having sufficiently large 
multiplicities are Cohen--Macaulay using Alexander duality. 
\end{abstract}
\maketitle
\section{Introduction}
\par 
Throughout this paper, let $S= k[X_1,\ldots,X_n]$
be a homogeneous polynomial ring over a field $k$ with $\deg X_i =1$. 
For a simplicial complex $\Delta$ on vertex set 
$[n]=\{1,\ldots,n\}$ (note that $\{i\} \in \Delta$ for all $i$), 
$k[\Delta] = k[X_1,\ldots,X_n]/I_{\Delta}$ is called 
the \textit{Stanley--Reisner ring} of $\Delta$, where $I_{\Delta}$ 
is an ideal generated by all square-free monomials 
$X_{i_1}\cdots X_{i_p}$ such that 
$\{i_1,\ldots,i_p\} \notin \Delta$. 
The ring $A = k[\Delta]$ is a homogeneous reduced ring with the 
unique homogeneous maximal ideal $\frm = (X_1,\ldots,X_n)k[\Delta]$ 
and the Krull dimension $d = \dim \Delta +1$. 
Let $e(A)$ denote the \textit{multiplicity} 
$e_0(\frm A_{\frm},A_{\frm})$ of $A$, 
which is equal to the number of facets (i.e., maximal faces)  
$F$ of $\Delta$ with $\dim F = d-1$. 
Also, we frequently call it the multiplicity of $\Delta$.  
Note that $\Delta$ is called \textit{pure} if all facets
of $\Delta$ have the same dimension. 
See \cite{BrHe, St} for more details. 
\par
Take a graded minimal free resolution
of a homogeneous $k$-algebra $A = S/I$ over $S:$
\[
 0 \to \bigoplus_{j \in \Z} S(-j)^{\beta_{p,j}(A)} 
\stackrel{\varphi_{p}}{\longrightarrow} \cdots 
\stackrel{\varphi_{2}}{\longrightarrow}  
\bigoplus_{j \in \Z} S(-j)^{\beta_{1,j}(A)}  
\stackrel{\varphi_{1}}{\longrightarrow}  
S \to A \to 0. 
\]
Then the \textit{initial degree} $\indeg A$ 
(resp. \textit{the relation type} $\rt(A)$) of $A$ 
is defined by 
$\indeg A = \min\{j \in \bbZ \,:\, \beta_{1,j}(A) \ne 0\}$ 
(resp. $\rt(A) = \max\{j \in \bbZ \,:\, \beta_{1,j}(A) \ne 0\}$). 
Also, $\reg A = \max\{j-i \in \bbZ \,:\, \beta_{i,j}(A) \ne 0\}$ is called  
the {\it Castelnuovo--Mumford regularity} of $A$. 
It is easy to see that $\reg A \ge \indeg A-1$, 
and $A$ has {\it linear resolution} if equality holds. 

\par
The main purpose of this paper is to prove the following theorems$:$ 

\par \vspace{2mm} \par \noindent
{\bf Theorem \ref{Main1}.}
Let $A = k[\Delta]$ be a Stanley--Reisner ring of Krull dimension $d \ge 2$. 
Put $\codim A = c$. 
If $e(A) \ge \sbinom{n}{c} - c$, then $A$ is Cohen--Macaulay. 

\par \vspace{2mm} \par \noindent
{\bf Theorem \ref{Main2}.}
Let $A = k[\Delta]$ be a Stanley--Reisner ring of Krull dimension $d \ge 2$. 
Put $\codim A = c$. 
Suppose that $\Delta$ is pure $($i.e., $A$ is equidimensional$)$. 
If $e(A) \ge \sbinom{n}{c} - 2c+1$, then $A$ is Cohen--Macaulay.

\par
It is easy to prove the above theorems in the case of 
$d=2$. 
When $d=2$, $A$ is Cohen--Macaulay if and only if $\Delta$ is connected. 
In fact, a disconnected graph has at 
most $\bbinom{n-1}{2}(=\bbinom{n}{2}-(n-2)-1)$ edges. 
This shows that Theorem \ref{Main1} is true in this case. 
Similarly, a disconnected graph without an isolated point 
has at most $\bbinom{n-2}{2}+1(=\bbinom{n}{2}-2(n-2))$ edges.
Indeed, such a graph is contained in a disjoint union of 
an $(n-i)$-complete graph and an $i$-complete graph for some $2 \le i \le n-2$. 
When $i=2$, the number of edges of the above union is just $\bbinom{n-2}{2}+1$. 
Thus we also get Theorem \ref{Main2} in this case. 

\par
The case $\indeg A = d$ and $c \ge 2$ is essential 
in the above two theorems. 
In order to prove Theorems \ref{Main1} and \ref{Main2} in this case, 
we consider their Alexander dual versions$:$

\par \vspace{2mm} \par \noindent 
{\bf Theorem \ref{AD-Main1}.}
Let $A = k[\Delta]$ be a Stanley--Reisner ring of Krull dimension $d \ge 2$. 
Suppose that $\indeg A = d$. 
If $e(A) \le d$, then $A$ has $d$-linear resolution. 
In particular, $\rt(A) = d$. 

\par \vspace{2mm} \par \noindent 
{\bf Theorem \ref{AD-Main2}.}
Let $A = k[\Delta]$ be a Stanley--Reisner ring of Krull dimension $d \ge 2$. 
Suppose that $\indeg A = \rt(A) = d$. 
If $e(A) \le 2d-1$, then $A$ has $d$-linear resolution. 
In particular, $a(A)< 0$.

\par \vspace{2mm} 
For a Stanley--Reisner ring $A$ with $\indeg A =\dim A =d$, 
it has $d$-linear resolution if and only if $a(R)< 0$. 
Thus the assertion of Theorem \ref{AD-Main2} could be seen as an analogy 
of the following: 
Let $R$ be a homogeneous integral domain over an algebraically closed 
field of characteristic $0$. 
If $e(R) \le 2\dim R-1$ and $\codim R \ge 2$, then $a(R)< 0$.

\par
In the last section, we will provide several examples related to 
the above results. 
\bigskip
\section{Complexes $\Delta$ with $e(k[\Delta]) \ge \bbinom{n}{c}-c$}

In this section, we use the following notation. 
Let $\Delta$ be a simplicial complex on $V=[n]$, and 
let $A = k[\Delta] = S/I_{\Delta}$ be the Stanley--Reisner ring of $\Delta$. 
Put $d = \dim A$, and $c = \codim A = n-d$. 
Note that $\bbinom{[n]}{d}$ denotes the family of 
all $d$-subsets of $[n]$. 

\par 
The main purpose of this section is to prove the following theorem. 

\begin{thm} \label{Main1}
Let $A = k[\Delta]$ be a Stanley--Reisner ring of Krull dimension $d \ge 2$. 
If $e(A) \ge \sbinom{n}{c} - c$, then $A$ is Cohen--Macaulay. 
\end{thm}

\par
Let us begin the proof of this theorem with the 
following lemmas. 

\begin{lemma} \label{Indeg}
If $e(A) \ge  \sbinom{n}{c} - c$, then $\indeg A \ge d$. 
\end{lemma}

\begin{proof}[\quad Proof]
Suppose that $\indeg A < d$. 
Take a squarefree monomial $M \in I_{\Delta}$ with $\deg M = d-1$. 
Then $(n-d+1)$ distinct squarefree monomials appear among $X_1 M,\ldots,X_n M$;
say $M_1,\ldots,M_{c+1}$. 
Let $F_i \in \bbinom{[n]}{d}$ corresponding to $M_i$, respectively. 
Then since no $F_i$ is contained in $\Delta$ we have 
\[
  e(A) \le \sbinom{n}{d} - (c+1). 
\]
This contradicts the assumption. 
\end{proof}

\begin{lemma} \label{IndegHigh}
Under the above notation, 
the following conditions are equivalent$:$
  \begin{enumerate}
   \item $\indeg A = d+1$.
   \item $e(A) = \sbinom{n}{d}$.
   \item $I_{\Delta} = (X_{i_1}\cdots X_{i_{d+1}}
\,:\,1 \le i_1 < \cdots < i_{d+1} \le n)$.
   \item $A$ has $(d+1)$-linear resolution. 
  \end{enumerate}
When this is the case, $A$ is Cohen--Macaulay with $\rt(A) = d+1$. 
\end{lemma}

\begin{proof}[\quad Proof]
See, e.g., \cite[Proposition 1.2]{TeYo}. 
\end{proof}

\par
Therefore we may assume 
that $\indeg A = d$ to prove Theorem \ref{Main1}. 

\begin{lemma} \label{Hyper}
Suppose $n = d+1$. 
If $e(A) \ge d$, then $A$ is a hypersurface.
\end{lemma}

\begin{proof}[\quad Proof]
Suppose that $A$ is not a hypersurface. 
Then we can write 
\[
I_{\Delta} = X_{i_1}\cdots X_{i_p}J
\]
for some monomial ideal $J(\ne R)$ with $\height J \ge 2$ 
since $\height I_{\Delta} =1$. 
In particular, $A$ is not Cohen--Macaulay. 
Thus $\indeg A \le d$ by Lemma \ref{IndegHigh}. 
Then $e(A) =p \le d-1$. 
This contradicts the assumption. 
\end{proof}

\par
Thus we may also assume that $c=\codim A \ge 2$. 
then let $\Delta^{*}$ be the {\it Alexander dual} of $\Delta$: 
\[
 \Delta^{*} = \{F \in 2^{V}\,:\, V \setminus F \notin \Delta\}. 
\]
Then $\Delta^{*}$ is a simplicial complex on the same vertex set $V$ 
of $\Delta$ for which the following properties are satisfied$:$ 

\begin{prop} \label{A-dualProp}
Under the above notation, we have 
\begin{enumerate}
 \item $\indeg k[\Delta^{*}]+\dim k[\Delta] = n$. 
 \item $\rt(k[\Delta^{*}]) = \bight I_{\Delta}$, 
where \[
 \bight I = \max\{\height \frp \,:\, \;
\text{$\frp$ is a minimal prime divisor of $I$}\}.
\]  
In particular, $\Delta$ is pure if and only if 
$\rt(k[\Delta^{*}]) = \indeg k[\Delta^{*}]$. 
 \item $\beta_{0,q^{*}}(I_{\Delta^{*}}) = e(k[\Delta])$, 
where $q^{*} = \indeg k[\Delta^{*}]$. 
 \item $(\Delta^{*})^{*} = \Delta$. 
\end{enumerate}
\end{prop}

\par
Also, the following theorem is fundamental. 
See \cite{EaRe} for more details. 

\begin{thm}[\textbf{Eagon--Reiner} \cite{EaRe}]
$k[\Delta]$ is Cohen--Macaulay if  and only if $k[\Delta^{*}]$ has 
linear resolution.
\end{thm}

\par
We want to reduce Theorem \ref{Main1} to 
its Alexander dual version. 
Let $\Delta^{*}$ be the Alexander dual of $\Delta$. 
Then $\indeg k[\Delta^{*}] = n-\dim k[\Delta]=c$ and 
$\dim k[\Delta^{*}] = n-\indeg k[\Delta] = n-d=c$. 
Also, since $\indeg k[\Delta^{*}]=\dim k[\Delta]=c$, 
we have 
\[
 e(k[\Delta^{*}]) = \sbinom{n}{c} - \beta_{0,c}(I_{\Delta^{*}})
 = \sbinom{n}{c} - e(A) \le c
\]
Therefore, it is enough to prove the following theorem. 

\begin{thm}[\textbf{Alexander dual version of Theorem \ref{Main1}}] 
\label{AD-Main1}
Let $A = k[\Delta]$ be a Stanley--Reisner ring of Krull dimension $d \ge 2$. 
Suppose that $\indeg A = d$. 
If $e(A) \le d$, then $A$ has $d$-linear resolution. 
In particular, $\rt(A) = d$. 
\end{thm}

\begin{proof}[\quad Proof]
(1) Put $a(A) = \sup\{p \in \bbZ \,:\, [H_{\frm}^d(A)]_p \ne 0\}$, the 
\textit{$a$-invariant} of $A$. 
From the assumption we obtain that
\[
a(A) +  d \le e(A) -1 \le d-1, 
\]
where the first inequality follows from e.g. \cite[Lemma 3.1]{HoMi}. 
Hence $a(A) < 0$. 
On the other hand, we have that 
$[H_{\frm}^{i}(A)]_j = 0$ for all $i$ and $j \ge 1$ since $A$ 
is a Stanley--Reisner ring. 
Then
\[
 \reg A = \inf\{p \in \bbZ \,:\, 
[H_{\frm}^i(A)]_j =0 \;\;\text{for all $i+j> p$}\}\le d-1=\indeg A-1. 
\]
This means that $A$ has $d$-linear resolution, as required. 
\end{proof}

\par
Now let us discuss a generalization of Theorem \ref{AD-Main1}. 
Let $A = S/I$ be an arbitrary homogeneous reduced $k$-algebra 
over a field $k$ of characteristic $p > 0$. 
The ring $A$ is called \textit{$F$-pure} if the Frobenius map 
$F \colon A \to A\;(a \mapsto a^p)$ is pure. 
It is known that a Stanley--Reisner ring is $F$-pure, and that 
if $A$ is $F$-pure then $[H_{\frm}^i(A)]_j =0$ for all $j \ge 1$. 
Thus the proof of Theorem \ref{AD-Main1}
involves that of the following proposition. 

\begin{prop} \label{F-pure-1}
Let $A = S/I$ be a homogeneous F-pure $k$-algebra. 
Put $\dim A =\indeg A =d \ge 2$.
If $e(A) \le d$, then $A$ has $d$-linear resolution. 
In particular, $\rt(A) =d$ and $a(A) < 0$. 
\end{prop}

\bigskip
\section{Complexes $\Delta$ with $e(k[\Delta]) \ge \bbinom{n}{c}-2c+1$}

\par
We use the same notation as in the previous section. 
For a face $G$ in $\Delta$ and $v \in V$, we put 
\begin{eqnarray*}
 \Delta_{V \setminus \{v\}} &=& \{F \in \Delta\,:\, v \notin F \}, \\
\Star_{\Delta} G  &=& \{F \in \Delta \,:\, F \cup G \in \Delta\}, \\
\link_{\Delta} G  &=& \{F \in \Delta \,:\, 
F \cup G \in \Delta, \, F \cap G = \emptyset \}. 
\end{eqnarray*}
\par
The main purpose of this section is to prove the following theorem. 

\begin{thm} \label{Main2}
Let $A = k[\Delta]$ be a Stanley--Reisner ring of Krull dimension $d \ge 2$. 
Put $c = \codim A$. 
Suppose that $\Delta$ is pure. 
If $e(A) \ge \sbinom{n}{c} - 2c+1$, then $A$ is Cohen--Macaulay. 
\end{thm}

\par
Now suppose that $c=1$ (resp. $\indeg A \ge d+1$). 
Then the assertion follows from Lemma \ref{Hyper} (resp. Lemma \ref{IndegHigh}).
Thus we may assume that $c \ge 2$ and $q=\indeg A \le d$. 
The following lemma corresponds to Lemma \ref{Indeg}. 

\begin{lemma} \label{Indeg2}
If $e(k[\Delta]) \ge \bbinom{n}{c}-2c+1$, then $\indeg k[\Delta] \ge d-1$, 
\par \noindent 
i.e., $(1)$ $\indeg k[\Delta] =d$ or $(2)$ $\indeg k[\Delta] =d-1$. 
\end{lemma}

\begin{proof}[\quad Proof]
Suppose that $\indeg k[\Delta] < d-1$. 
Take a squarefree monomial $M \in I_{\Delta}$ with $\deg M  = d-2$. 
Then there are $\bbinom{n-d+2}{2}$ squarefree monomials in degree $d$ 
in $I_{\Delta}$. 
Note $\bbinom{n-d+2}{2} = \bbinom{c+2}{2} \ge 2c$. 
Hence 
\[
  e(k[\Delta]) \le \bbinom{n}{c}-2c. 
\]
This contradicts the assumption. 
\end{proof}

\par
First, we consider the Alexander dual version of Theorem \ref{Main2} in the 
case of $\indeg k[\Delta]=d$. 
Namely, we will prove the following theorem. 

\begin{thm}[\textbf{Alexander dual version of Theorem \ref{Main2}, Case $(1)$}] 
\label{AD-Main2}
Let $A = k[\Delta]$ be a Stanley--Reisner ring of Krull dimension $d \ge 2$. 
Suppose that $\indeg A = \rt(A) = d$. 
If $e(A) \le 2d-1$, then $A$ has $d$-linear resolution. 
In particular, $a(A)< 0$. 
\end{thm}

\par
The proof of the above theorem can be reduced to that of the following theorem, 
which is a key result in this paper. 

\begin{thm} \label{Key}
Let $A = k[\Delta]$ be a Stanley--Reisner ring of Krull dimension $d \ge 2$. 
Suppose that $\rt(A) \le d$. 
If $e(A) \le 2d-1$, then $\reg A \le d-1$, equivalently, 
$\widetilde{H}_{d-1}(\Delta) = 0$. 
\end{thm}

\begin{proof}[\quad Proof]
Put $e = e(A)$. 
Let $\Delta'$ be the subcomplex that is spanned by all facets of dimension $d-1$.
Replacing $\Delta$ with $\Delta'$, we may assume that $\Delta$ is pure. 
\par
We use induction on $d = \dim A \ge 2$. 
First suppose $d=2$. 
The assumption shows that $\Delta$ does not contain the boundary complex 
of a triangle. 
Hence $\widetilde{H}_1(\Delta) = 0$
since $e(A) \le 3$. 
\par
Next suppose that $d \ge 3$, and that the assertion holds for any complex 
the dimension  of which is less than $d-1$. 
Assume that $\widetilde{H}_{d-1}(\Delta) \ne 0$. 
Take one $\Delta$ whose multiplicity is minimal
among the multiplicities of those complexes. 
Then $\Delta$ does not contain any {\it free face} (see \cite{Hud}).  
That is, every face that is not a facet is contained in at least two facets. 
Indeed, suppose that $\Delta$ contains a free face (say, $G$) and 
put $\Delta' = \Delta \setminus \{F \in \Delta \,:\, F \supseteq G\}$. 
Then since $G$ is a free face of $\Delta$, $\Delta'$ is homotopy equivalent to 
$\Delta$ and $e(k[\Delta'])=e(k[\Delta])-1$. 
In particular, $\widetilde{H}_{d-1}(\Delta')\cong \widetilde{H}_{d-1}(\Delta)\ne 0$.
This contradicts the minimality of $e(k[\Delta])$. 
\par
First consider the case of $\rt(A) =d$. 
Take a generator $X_{i_1}\cdots X_{i_d}$ of $I_{\Delta}$. 
For every $j=1,\ldots,d$, each $G_j = \{i_1,\ldots,\widehat{i_j},\ldots,i_d\}$ 
is contained in at least two facets as mentioned above. 
Then $e(A) \ge 2d$ since those facets are 
different from each other.  
This is a contradiction. 
\par
Next we consider the case of $\rt(A) < d$. 
Take a Mayer--Vietoris sequence with respect to $\Delta = \Delta_{V \setminus \{n\}}
\cup \Star_{\Delta}\{n\}$ as follows$:$
\[
\widetilde{H}_{d-1}(\Delta_{V \setminus \{n\}}) 
\oplus 
\widetilde{H}_{d-1}(\Star_{\Delta}\{n\})
\longrightarrow 
\widetilde{H}_{d-1}(\Delta)
\longrightarrow 
\widetilde{H}_{d-2}(\link_{\Delta}\{n\}). 
\]
The minimality of $e(k[\Delta_{V \setminus \{n\}}])$ yields 
that $\widetilde{H}_{d-1}(\Delta_{V \setminus \{n\}}) =0$ 
since $e(k[\Delta_{V \setminus \{n\}}]) < e(k[\Delta])$. 
On the other hand, it is known that 
$\widetilde{H}_{i}(\Star_{\Delta}\{n\}) =0$ for all $i$. 
Hence $\widetilde{H}_{d-1}(\Delta) \hookrightarrow  
\widetilde{H}_{d-2}(\link_{\Delta}\{n\})$. 
In particular, $\widetilde{H}_{d-2}(\link_{\Delta}\{n\}) \ne 0$. 
\par
Set $\Delta' = \link_{\Delta}\{n\}$. 
Then $\Delta'$ is a complex on $V \setminus \{n\}$ 
such that $\dim k[\Delta']= d-1$ and
$\rt(k[\Delta']) \le \rt(k[\Delta]) \le d-1$. 
In order to apply the induction hypothesis to $\Delta'$, 
we want to see that $e(k[\Delta']) \le 2d-3$. 
In order to do that, we consider $e(k[\Delta_{V \setminus \{n\}}])$. 
As $\Delta \ne \Star_{\Delta}\{n\}$, one can take 
$F = \{i_1,\cdots, i_p,n\} \notin \Delta$ for some $p \le d-2$ such that 
$X_{i_1}\cdots X_{i_p}X_n$ is a generator of $I_{\Delta}$. 
Then $G := \{i_1,\ldots,i_p\} \in \Delta$, but it is not a facet of $\Delta$. 
Thus it is contained in at least 
two facets of $\Delta$, each of which does not contain $n$. 
Hence $e(k[\Delta_{V \setminus \{n\}}]) \ge 2$. 
Thus we get 
\[
e(k[\Delta']) 
 =  e(k[\Star_{\Delta}\{n\}]) 
 =  e(k[\Delta])- e(k[\Delta_{V \setminus \{n\}}]) 
\le 2d-3. 
\]
By induction hypothesis, we have 
$\widetilde{H}_{d-2}(\link_{\Delta}\{n\})=0$.
This is a contradiction. 
\end{proof}

\par
Next, we consider the Alexander dual version of Theorem \ref{Main2} in the case 
of $\indeg k[\Delta]=d-1$. 
Namely, we must prove the following proposition. 

\begin{prop}[\textbf{Alexander dual version of Theorem \ref{Main2}, Case $(2)$}] 
\label{Omake}
Let $A = k[\Delta]$ be a Stanley--Reisner ring of Krull dimension $d \ge 2$. 
Suppose that $\indeg A = \rt(A) = d-1$. 
If $\mu(I_{\Delta}) \ge \sbinom{n}{d-1} -2d+3$, then 
$A$ has $(d-1)$-linear resolution with $e(A) =1$. 
\end{prop}

\begin{proof}[\quad Proof]
First we show that $e(A)=1$. 
Now suppose that $e(A)\ge 2$. 
Then there exist at least two facets $F_1$ and $F_2$ with $\#(F_1)=\#F(F_2)=d$. 
This implies that $f_{d-2}(\Delta) \ge 2d-1$. 
However, by the assumption, we have 
\[
 f_{d-2}(\Delta) = \bbinom{n}{d-1}-\beta_{0,d-1}(I_{\Delta}) = 
 \bbinom{n}{d-1}-\mu(I_{\Delta}) \le 2d-3. 
\]
This is a contradiction. 
Hence we get $e(A)=1$. 
\par
In order to prove that $A$ has $(d-1)$-linear resolution, 
it is enough to show that 
$\beta_{i,j}(A)=0$ for all $i \ge c$ and $j \ge i+d-1$
by \cite[Theorem 5.2]{Sch1}. 
Also, it suffices to show that 
$\widetilde{H}_{d-1}(\Delta) = \widetilde{H}_{d-2}(\Delta)=
\widetilde{H}_{d-2}(\Delta_W)=0$
for all subsets $W \subset V$ with $\#(W)=n-1$
by virtue of Hochster's formula on the Betti numbers:
\[
 \beta_{i,j}(A) = \sum_{\begin{subarray}{c} W \subseteq V \\
\#(W) = j\end{subarray}} 
\!\!\!\! \dim_k \widetilde{H}_{j-i-1}(\Delta_W;k). 
\]
\par \vspace{2mm}
\begin{tabular}{rl}
{\bf Claim 1}. &  $\widetilde{H}_{d-1}(\Delta) = \widetilde{H}_{d-2}(\Delta)=0$. 
\end{tabular}
\par \vspace{2mm}
Since $\rt(A) \le d-1 \le d$ and $e(A) =1 \le 2d-1$, we have 
$\widetilde{H}_{d-1}(\Delta)=0$ by Theorem \ref{Key}. 
Now let $F=\{1,2,\ldots,d\}$ be the unique facet with $\#(F)=d$. 
Consider a simplicial subcomplex $\Delta' := \Delta \setminus \{F,G\}$ 
where $G=\{1,2,\ldots,d-1\}$. 
Then $\dim k[\Delta']=d-1$ and $e(k[\Delta']) \le 2d-4 \le 2(d-1)-1$. 
Also, since $\rt(k[\Delta']) \le \rt(k[\Delta])\le d-1$, applying 
Theorem \ref{Key} to $\Delta'$, we obtain that 
$\widetilde{H}_{d-2}(\Delta) \cong \widetilde{H}_{d-2}(\Delta')=0$, as required. 
\par \vspace{2mm}
\begin{tabular}{rl}
{\bf Claim 2}. & 
$\widetilde{H}_{d-2}(\Delta_W)=0$
for all subsets $W \subset V$ with $\#(W)=n-1$. 
\end{tabular}
\par \vspace{2mm}
Let $W$ be a subset of $V$ such that $\#(W)=n-1$. 
Put $\{a\} = V \setminus W$. 
If $a$ is not contained in $F$, then $\widetilde{H}_{d-2}(\Delta_W)=0$ by the
similar argument as in the proof of the previous claim. 
So we may assume that $a \in F$. 
Then $\dim k[\Delta_W]=d-1$ and $e(k[\Delta_W]) \le (d-3)+1=d-2 \le 2(d-1)-1$. 
Also, since $\rt(k[\Delta_W]) \le d-1$, we have 
$\widetilde{H}_{d-2}(\Delta_W) =0$ by Theorem \ref{Key} again.  
\par \vspace{2mm} \par \noindent
Hence $k[\Delta]$ has $(d-1)$-linear resolution, as required. 
\end{proof}

\begin{exam}  \label{OmakeEx}
Let $\rho$, $d$ be an integers with $0 \le \rho \le d-3$. 
Let $\Delta$ be a simplicial complex on $V=[n]$
spanned by $F=\{1,2,\ldots,d\}$, any distinct $\rho$ elements 
from $\bbinom{[n]}{d-1} \setminus \bbinom{[d]}{d-1}$ and 
all elements of $\bbinom{[n]}{d-2}$. 
Then $\dim k[\Delta]=d$, $\indeg k[\Delta]=\rt(k[\Delta])=d-1$. 
Also, we have 
\[
 \mu(I_{\Delta}) = \beta_{0,d-1}(I_{\Delta}) = \bbinom{n}{d-1}-\rho-d \ge 
\bbinom{n}{d-1}-2d+3.
\]
Hence $\Delta$ satisfies the assumption of the above proposition. 
\end{exam}

\par \vspace{3mm}
On the other hand, we have no results for F-pure $k$-algebras 
corresponding to Theorem \ref{AD-Main2}. 
But we remark the following. 

\begin{remark} \label{Domain-remark}
As mentioned in the introduction, if $A$ is a homogeneous integral domain 
over an algebraically closed field of $\chara k =0$ with $\codim A \ge 2$ and
$e(A) \le 2d-1$ then one has $a(A) < 0$. 
In fact, it is known that an inequality 
\[
  a(A) + d \le \left\lceil \frac{e(A)-1}{\codim A} \right\rceil
\]
holds; see e.g., the remark after Theorem $3.2$ in \cite{HoMi}. 
Moreover, Professor Chikashi Miyazaki told us that this inequality is 
also true in positive characteristic. 
\end{remark}

\begin{quest} \label{Fpure2}
Let $A = k[A_1]$ be a homogeneous F-pure, equidimensional $k$-algebra. 
Put $\dim A = \indeg A =d \ge 2$. 
If $e(A) \le 2d-1$, then does $a(A)<0$ hold? 
\end{quest}

\bigskip 
\section{Buchsbaumness}

\par
A Stanley--Reisner ring $A = k[\Delta]$ is \textit{Buchsbaum} if and only if 
$\Delta$ is pure and $k[\link_{\Delta}\{i\}]$ is Cohen--Macaulay for every 
$i \in [n]$. 
As an application of Theorem \ref{Main2}, we can provide sufficient 
conditions for $k[\Delta]$ to be Buchsbaum.

\begin{prop} \label{Bbm}
Let $A = k[\Delta]$ be a Stanley--Reisner ring of Krull dimension 
$d \ge 3$. 
Suppose that $\Delta$ is pure, $\indeg A = d$ and 
$e(A) \ge \bbinom{n}{c}-2c$. 
Then 
\begin{enumerate}
 \item $e(k[\link_{\Delta}\{i\}]) \ge \bbinom{n-1}{c}-2c$ for all $i$.  
 \item If $\height [I_{\Delta}]_dS  \ge 2$, then $A$ is Buchsbaum. 
 \item If $\rt(A) = d$, then $A$ is Buchsbaum. 
\end{enumerate}
\end{prop}

\begin{proof}[\quad Proof]
We may assume that $c \ge 2$, $e(A) = \bbinom{n}{c}-2c$, 
and that 
$\Delta \ne \Star_{\Delta}\{i\}$ for every $i \in [n]$. 
Put $\Gamma_i = \link_{\Delta}\{i\}$ for each $i \in [n]$. 
\par
(1)  
We first show the following claim. 
\par \vspace{1mm}
\par \noindent 
{\bf Claim:} $e(A) \le \bbinom{n}{d}- 
\left\{\bbinom{n-1}{d-1}-e(k[\Gamma_i])\right\}$ for every $i \in [n]$. 
Also, equality holds if and only if $i \in F$ holds for all $F \in \bbinom{[n]}{d}
\setminus \Delta$. 

\par \vspace{1mm}
Put $W_i = \left\{F \in \bbinom{[n]}{d}\,:\, i \in F \notin
\Star_{\Delta}\{i\} \right\}$. 
Then $\#(W_i) \le \# (\bigcup_{i=1}^n W_i)$ implies that 
\[
 \bbinom{n-1}{d-1} - e(k[\Gamma_j]) \le \bbinom{n}{d}-e(A),
\]
as required. 
Also, equality holds if and only if $W_i = \bigcup_{i=1}^n W_i$, 
that is,  $i \in F$ holds for all $F \in \bbinom{[n]}{d}
\setminus \Delta$.
\par
Now suppose that $e(k[\Gamma_i]) \le \bbinom{n-1}{d-1} - 2c-1$
for some $i \in [n]$. 
Then the claim implies that $e(A) \le \bbinom{n}{d}-2c-1$, which 
contradicts the assumption. 
Thus we get (1). 
\par 
(2) Suppose that $\height [I_{\Delta}]_dS \ge 2$. 
Then there is no element $i \in [n]$ for which $i \in F$ holds for all 
$F \in \bbinom{[n]}{d} \setminus \Delta$. 
Thus the claim yields that 
\[
 \bbinom{n}{d}-2c = e(A) \le \bbinom{n}{d}- 
\left[\bbinom{n-1}{d-1} - e(k[\Gamma_i]) \right]-1, 
\]
that is, $e(k[\Gamma_i]) \ge \bbinom{n-1}{d-1}-2c+1$ for every $i \in [n]$. 
Also, we note that $\Gamma_i$ is pure and 
$\indeg k[\Gamma_i] = \dim  k[\Gamma_i] =d-1$. 
Applying Theorem \ref{Main2} to $k[\Gamma_i]$, we obtain that 
$k[\Gamma_i]$ is Cohen--Macaulay. 
Therefore $A$ is Buchsbaum since $\Delta$ is pure. 
\par
(3) Now suppose that $A$ is \textit{not} Buchsbaum. 
Then since $\height [I_{\Delta}]_d S =1$, one can take $i \in [n]$
for which $i \in F$ holds for all $F \in \bbinom{[n]}{d} \setminus \Delta$. 
We may assume $i=n$. 
Then $\{1,\ldots, \widehat{i},\ldots,d+1\} \in \Delta$ for all 
$i \in [d+1]$ because $n-1 \ge d+1$. 
This means that $X_1\cdots X_{d+1}$ is a generator of $I_{\Delta}$; thus 
$\rt(A) = d+1$. 
\end{proof}

\bigskip
\section{Examples}

Throughout this section, let $c$, $d$ be given integers with $c,\,d \ge 2$. 
Set $n=c+d$.  

\begin{exam} \label{Thm-sample}
Put $F_{i,j} = \{1,2,\ldots,\widehat{i},\ldots, d,j\}$ for each 
$i=1,\ldots,d; j = d+1,\ldots,n$. 
For a given integers $e$ with $1 \le e \le cd$, we choose $e$ faces 
(say, $F_1,\ldots, F_e$) from 
$\{F_{i,j}\,:\, 1 \le i \le d,\,d+1 \le j \le n\}$, 
which is a simplicial join of $2^{[d]} \setminus \{[d]\}$ and $c$ points. 
\par
Let $\Delta$ be a simplicial complex spanned by $F_1,\ldots,F_e$ and 
all elements of $\bbinom{[n]}{d-1}$. 
Then $k[\Delta]$ is a $d$-dimensional Stanley--Reisner ring with 
$\indeg k[\Delta] = \rt(k[\Delta])=d$ and $e(k[\Delta]) = e$. 
\par
In particular, when $e \le 2d-1$, 
$k[\Delta]$ has $d$-linear resolution by Theorem \ref{AD-Main2}. 
Thus 
Also, the Alexander dual complexes of them provide examples 
satisfying hypothesis of Theorem \ref{Main2}.
\end{exam}

\par
The following example shows that 
the assumption \lq\lq $e(A) \le 2d-1$'' is 
optimal in Theorem \ref{AD-Main2}. 

\begin{exam} \label{Exam-notlin}
There exists a complex $\Delta$ on $V = [n]$ ($n=d+2$) 
for which $k[\Delta]$ does not have 
$d$-linear resolution with $\dim k[\Delta] = \indeg k[\Delta] 
= \rt(k[\Delta]) = d$ and $e(k[\Delta]) = 2d$. 
\par
In fact, put $n=d+2$. 
Let $\Delta_0$ be a complex on $V = [n]$ such that 
$k[\Delta_0]$ is a complete intersection 
defined by $(X_1\cdots X_d,X_{d+1}X_{d+2})$. 
Also, let $\Delta$ be a complex on $V$ that is spanned 
by all facets of $\Delta_0$ and all elements of $\bbinom{[n]}{d-1}$ $:$
\[
 I_{\Delta} = 
 (X_1\cdots X_d)S+ 
 (X_{i_1}\cdots X_{i_{d-2}}X_{d+1}X_{d+2}\,:\, 
 1 \le i_1 < \cdots < i_{d-2} \le d )S. 
\]
Then $\widetilde{H}_{d-1}(k[\Delta]) \cong \widetilde{H}_{d-1}(k[\Delta_0]) \ne 0$ 
since $a(k[\Delta_0]) = 0$. 
Hence $k[\Delta]$ does not have linear resolution. 
\end{exam}

\begin{remark}
The above example is obtained 
by considering the case $c=2,\,e=2d$ in Example $\ref{Thm-sample}$. 
\end{remark}

\par
The next example shows that 
the assumption \lq\lq $\rt(A) =d$'' is 
\textit{not} superfluous in Theorem \ref{AD-Main2}. 

\begin{exam} \label{Exam-rt}
Suppose that $d+1 \le e \le \bbinom{n}{d}-1$. 
There exists a simplicial complex $\Delta$ on $V = [n]$ such that 
$\dim k[\Delta] = \indeg k[\Delta]=d$, $\rt(k[\Delta]) = d+1$ and $e(k[\Delta])=e$. 
In particular, $k[\Delta]$ does not have $d$-linear resolution. 
\par
In fact, put $\mathcal{F}  =\bbinom{[n]}{d}\setminus \bbinom{[d+1]}{d}$. 
Let $\Delta_0$ be a simplicial complex on $V$ such that 
\[
I_{\Delta_0} = (X_1\cdots X_{d}X_{d+1})S+
(X_{i_1}\cdots X_{i_d}\,:\, \{i_1,\ldots,i_d\} \in 
\mathcal{F})S.
\]
Then $\dim k[\Delta_0] = \indeg k[\Delta_0]=d$, $\rt(k[\Delta_0]) =d+1$, and 
$e(k[\Delta_0])=d+1$. 
\par
For a given integer $e$ which satisfies above condition, 
one obtains the required simplicial complex by adding any $(e-d-1)$ distinct 
$d$-subsets of $2^{[n]}$ that is not contained in $\bbinom{[d+1]}{d}$ to $\Delta_0$. 
\end{exam}

\begin{remark} \label{TuranNum}
Now let $\Delta$ be a simplicial complex on $V=[n]$. 
Set $A = k[\Delta]$. 
Suppose that $\dim A=\indeg A=d \ge 2$. 
Then one can easily see that $d \le \rt(A) \le d+1$; 
$\rt(A) = d$ $($resp. $d+1$$)$ if $1 \le e(A) \le d$ 
$($resp. $e(A) = \bbinom{n}{d}$ $)$. 
So we put 
\[
f(n,d) = 
\min\left\{m \in \bbZ\,:\, 
\begin{array}{l}
\rt k[\Delta]=d+1 \;\;\text{for all $(d-1)$-dimensional} \\
\text{complexes $\Delta$ on $V$ with $\indeg k[\Delta] = d$} \\
\text{and $e(k[\Delta]) \ge m$}
\end{array}
\right\}
\]
Then $f(n,d) \ge cd+1$ by Example $\ref{Thm-sample}$. 
From the definition of $f(n,d)$, one can easily see that 
there exists a simplicial complex $\Delta$ which satisfies $\rt(k[\Delta]) = d$ and 
$e(k[\Delta])=e$ for each $e$ with $d+1 \le e \le f(n,d)-1$. 
On the other hand, by virtue of Example $\ref{Exam-rt}$, 
one can also find a simplicial complex $\Delta$ which satisfies 
$\rt(k[\Delta]) = d+1$ and
$e(k[\Delta])=e$ for each $e$ with $d+1 \le e \le \bbinom{n}{d}-1$. 
\par
It seems to be difficult to determine $f(n,d)$ in general. 
Let $T(n,p,k)$ be the so-called \textit{Turan number}. Then we have 
\[
 f(n,d) = \bbinom{n}{d}- T(n,d+1,d). 
\]
In particular, we get 
\begin{equation} \label{Turan}
 f(n,2) = \left\{
\begin{array}{ll}
\frac{n^2}{4}+1, & \;\text{if $n$ is even$;$} \\[2mm]
\frac{n^2-1}{4}+1, & \;\text{otherwise} 
\end{array}
\right. 
\end{equation}
by Turan's theorem $($e.g., \cite[Theorem 7.1.1]{Di}$)$
However, no formula is known for $T(n,4,3)$; see \cite[pp.1320]{Fr}. 
\end{remark}

\par
In the rest of this section, we show that 
the purity of $\Delta$ is very strong condition in Theorem \ref{AD-Main2}. 

\begin{prop} \label{Pure}
Then the following conditions are equivalent$:$
\begin{enumerate}
 \item There exists a $d$-dimensional Stanley--Reisner ring $k[\Delta]$ 
such that $\Delta$ is pure, $\indeg k[\Delta] = d$ 
and $e(k[\Delta])=e \le 2d-1$. 
 \item $n=d+2$, $d \le 5$ and $(d,e)$ is one of the following pairs$:$
\[
(2,2),\,(2,3),\,(3,4),\,(3,5),\,(4,6),\,(4,7),\,(5,9).
\]
\end{enumerate}
\end{prop}

\par
To prove the proposition, we need the following lemma. 

\begin{lemma} \label{PureVertex}
Let $A = k[\Delta]$ be a $d$-dimensional Stanley--Reisner ring which is not a 
hypersurface. 
Suppose that $\Delta$ is pure and $\indeg A =d \ge 3$. 
Then there exists a vertex $i \in [n]$ such that 
$e(k[\Delta_{V \setminus \{i\}}]) \ge 2$. 
\end{lemma}

\begin{proof}[\quad Proof]
Note that $n \ge d+2$ by the assumption. 
Put $e = e(A)$.
Suppose that $e(k[\Delta_{V \setminus \{i\}}]) =1$ for all $i$. 
Then since there exist $(e-1)$ facets containing $i$ for each $i \in [n]$, 
we have 
\[
  (d+2)(e-1)\le n(e-1) \le de;
\]
hence $e \le \frac{d+2}{2}$. 
\par
On the other hand, by counting the number of subfacets (i.e., the maximal faces 
among all  faces except facets) of $\Delta$ we get 
\[ 
  de \ge \bbinom{n}{d-1}
\]
since $\indeg A =d$ and $\Delta$ is pure. 
It follows from these inequalities that 
\[
 \frac{d(d+2)}{2} \ge de \ge \bbinom{n}{d-1} \ge \bbinom{d+2}{d-1} = \bbinom{d+2}{3}. 
\]
Hence $d \le 2$. 
This is a contradiction. 
\end{proof}

\begin{proof}[\quad Proof of Proposition $\ref{Pure}$]
We first show $(1) \Longrightarrow (2)$. 
Let $A = k[\Delta]$ be a $d$-dimensional Stanley--Reisner ring for which 
$\Delta$ is pure, $\indeg A = d$, and $e = e(A) \le 2d-1$. 
We may assume that $d \ge 3$.
Since $\Delta$ is pure, any subfacet is contained in some $d$-subset of $\Delta$. 
By counting the number of subfacets that contain $n$, we obtain that 
\[
 \bbinom{n-1}{d-2} \le \left(e-e(k[\Delta_{V \setminus \{n\}}]\right)(d-1)
\le (e-2)(d-1), 
\]
where the last inequality follows from Lemma \ref{PureVertex}. 
\par 
Now let us see that $n=d+2$. 
Suppose that $n \ge d+3$. 
Then we get 
\[
 \bbinom{d+2}{4} \le \bbinom{n-1}{d-2} \le (e-2)(d-1) \le (2d-3)(d-1)
\]
by the assumption. 
This implies that $d \le 4$. 
\par 
First we consider the case of $d=4$. 
Then $n=d+3=7$, $e=2d-1 =7$. 
Let $\{F_1,\ldots,F_7\}$ be the set of facets of $\Delta$.  
Since $e(k[\Delta_{V \setminus \{7\}}])=2$, 
we may assume that $7 \in F$ if and only if $1 \le i \le 5$. 
Note that $F_i$ contains only one subfacet that does not contain $7$
for each $1 \le i \le 5$. 
On the other hand, one can find at most $4\times 2$ subfacets 
as faces of $F_6$ or $F_7$.
Therefore the total number of subfacets that do not contain $7$ is at most $13$. 
However the number of all subfacets which do not contain $7$ 
is $\bbinom{7-1}{4-1}=20$ since $\indeg A = 4$. 
This is a contradiction. 
\par
By the similar observation as in the case of $d=4$, 
one can prove that the case of $d=3$ does not occur. 
Therefore we conclude that $n=d+2$. 
\par
Under the assumption that $n=d+2$, let us determine $(d,e)$. 
Let $\Delta^{*}$ be the Alexander dual of $\Delta$ and put $R=k[\Delta^{*}]$. 
Then $R$ is a two-dimensional Stanley--Reisner ring with $\indeg R = 2$. 
Also, $\rt(R) = \indeg R=2$ since $\Delta$ is pure. 
Thus by virtue of Turan's theorem (see Eq. \ref{Turan}), 
we have 
\[
 \bbinom{d+2}{2}-e = e(R) \le f(d+2,2)-1 = \left\lfloor\frac{(d+2)^2}{4}\right\rfloor, 
\]
where $\lfloor a \rfloor$ denotes the maximum integer that does not exceed $a$. 
Namely, we have 
\[
 2d-1 \ge e \ge \left\lfloor \frac{(d+1)^2}{4} \right\rfloor. 
\]
It immediately follows from here that $(d,e)$ is one of the pairs 
listed above. 
\par
Conversely, in order to prove $(2) \Longrightarrow (1)$, it is enough to 
find $(n,e')$-graphs (i.e., $1$-dimensional simplicial complexes $\Gamma$ 
on $[n]$ with $e'$ edges) which does not contain any triangle 
for each $(n,e') =(4,4)$, $(4,3)$, $(5,6)$, $(5,5)$, $(6,9)$, $(6,8)$, $(7,12)$. 
Those complexes will be given in the following example. 
\end{proof}

\begin{exam} \label{PureDual}
There exists a $1$-dimensional simplicial connected complex  
$\Gamma$ on $[n]$ with  with $e(k[\Gamma])=e'$ and $\rt(k[\Gamma])=2$ 
for each $(n,e') =(4,4)$, $(4,3)$, $(5,6)$, $(5,5)$, $(6,9)$, $(6,8)$, 
$(7,12)$.  
Put 
\begin{eqnarray*}
S_{4,4} &=& \{[12],[14],[23],[34]\}, \\
S_{4,3} &=& \{[12],[23],[34] \}, \\
S_{5,6} &=& \{[12],[14],[23],[25],[34],[45]\},\\
S_{5,5} &=& \{[12],[14],[23],[34],[45]\}, \\
S_{6,9} &=& \{[14],[15],[16],[24],[25],[26],[34],[35],[36]\}, \\
S_{6,8} &=& \{[12],[14],[23],[25],[34],[36],[45],[56]\}, \\
S_{7,12} &=& \{[15],[16],[17],[25],[26],[27],[35],[36],[37],[45],[46],[47]]\},
\end{eqnarray*}
where $[i_1i_2\cdots i_p]$ means $\{i_1,i_2,\ldots,i_p\}$. 
\par 
Let $\Gamma_{n,e'}$ be a simplicial complex spanned by $S_{n,e'}$. 
Then $k[\Gamma_{n,e'}]$ is a two-dimensional Cohen--Macaulay Stanley--Reisner ring 
with $e(k[\Gamma])=e'$ and $\rt(k[\Gamma])=2$. 
Note that when $e' = f(n,2)-1$, $\Gamma_{n,e'}$ is the so-called 
\textit{Turan graph} $T^2(n)$, that is, it is the unique complete bipartite graph 
on $[n]$ whose two partition sets differ in size by at most $1$. 

\par
Let $\Delta_{d,e}$ be the Alexander dual complex of $\Gamma_{n,e'}$  
where $d=n-2$ and $e=\bbinom{d+2}{2}-e'$. 
Namely, $\Delta_{d,e}$ is the complex spanned by $T_{d,e}$, respectively$:$
\begin{eqnarray*}
T_{2,2} &=& \{[13],[24]\}, \\
T_{2,3} & =& \{[13],[23],[24]\}, \\
T_{3,4} &=& \{[124],[135],[234],[245]\}, \\
T_{3,5} &=&  \{[124],[134],[135],[234],[245]\}, \\
T_{4,6} &=& \{[1234],[2345],[3456],[4561],[5612],[6123]\}, \\
T_{4,7} &=& \{[1235],[1246],[1345],[1356],[2345],[2346],[2456]\}, \\
T_{5,9} &=&  \{[12345],[12346],[12347],[12567],[13567],[14567],[23567], \\
 & &[24567],[34567]\}. 
\end{eqnarray*}
Then $A = k[\Delta_{d,e}]$ is a $d$-dimensional equidimensional 
Stanley--Reisner ring with $\indeg A = d$ and $e(A) = e$ for 
each $(d,e) = (2,2)$, $(2,3)$, $(3,4)$, $(3,5)$, $(4,6)$, $(4,7)$, $(5,9)$.  
\end{exam}

\begin{cor} \label{Bbm-Pure}
Let $A = k[\Delta]$ be a $d$-dimensional Buchsbaum Stanley--Reisner ring which 
is not a hypersurface. 
Suppose that $\indeg A =d \ge 3$. 
Then $d=3$ and $\Delta$ is isomorphic to 
a simplicial complex spanned by $\{[124],[134],[135],[235],[245]\}$. 
\end{cor}

\begin{proof}[\quad Proof]
Since $A$ is Buchsbaum and $\indeg A =d$ we have 
\[
e=e(A) \ge \frac{c+d}{d}\bbinom{c+d-2}{d-2}
\]
by \cite[Proposition 2.1]{TeYo}. 
Also, $n=d+2$ by Proposition \ref{Pure} since $\Delta$ is pure. 
Thus
\[
2d-1 \ge e \ge \frac{d+2}{d}\bbinom{d}{d-2} = \frac{(d+2)(d-1)}{2}. 
\]
This implies that $d \le 3$, and thus $d=3$ and $e=5$. 
Then one can easily see that $\Delta$ is isomorphic to the 
complex spanned by $\{[124],[134],[135],[235],[245]\}$, 
which is the Alexander dual complex of a $1$-dimensional connected complex 
spanned by $\{[12],[23],[34],[45]\}$. 
\end{proof}

\end{document}